\documentclass{amsart}

\usepackage{amssymb}
\usepackage{graphicx}

\title{Inhomogeneous {S}trichartz estimates}
\author{Damiano Foschi}
\address{Dipartimento di Matematica Pura e Applicata\\
  Universit\`a di L'Aquila}
\email{foschi@univaq.it}
\urladdr{http://univaq.it/\textasciitilde foschi/}
\date{December 8, 2003}


\newcommand{\epsi}{\varepsilon}

\newcommand{\Z}{\mathbb{Z}}  
\newcommand{\R}{\mathbb{R}}  
\newcommand{\iC}{\mathrm{i}} 
\newcommand{\Eu}{\mathrm{e}} 

\newcommand{\les}{\lesssim}
\newcommand{\ges}{\gtrsim}

\newcommand{\sss}{\approx}

\newcommand{\abs}[1]{\left|#1\right|}
\newcommand{\norm}[1]{\left\|#1\right\|}
\newcommand{\tonde}[1]{\left(#1\right)}
\newcommand{\quadre}[1]{\left[#1\right]}
\newcommand{\graffe}[1]{\left\{#1\right\}}

\newcommand{\TOnde}[1]{\Bigl(#1\Bigr)}
\newcommand{\GRaffe}[1]{\Bigl\{#1\Bigr\}}
\newcommand{\NOrm}[1]{\Bigl\|#1\Bigr\|}

\newcommand{\ip}[1]{\langle #1 \rangle} 

\newenvironment{sistema}%
{\left\{\begin{aligned}}{\end{aligned}\right.}

\newcommand{\de}{\partial}
\renewcommand{\d}{\,{\rm d}}

\newcommand{\union}{\cup}        
\newcommand{\intersection}{\cap} 

\newcommand{\1}[1]{\frac{1}{#1}}

\newcommand{\mB}{{\mathcal B}}
\newcommand{\mE}{{\mathcal E}}
\newcommand{\mF}{{\mathcal F}}
\newcommand{\mQ}{{\mathcal Q}}
\newcommand{\mX}{{\mathcal X}}
\newcommand{\pt}{{\widetilde{p}}}
\newcommand{\qt}{{\widetilde{q}}}
\newcommand{\rt}{{\widetilde{r}}}
\newcommand{\Qt}{{\widetilde{Q}}}
\newcommand{\Rt}{{\widetilde{R}}}

\newcommand{\LtLX}[3][\R]{{L_t^{#2} (#1; L_X^{#3})}}

\newcommand{\Elocal}{{\mE_{\rm local}}}

\newcommand{\ie}{{i.e.\ }}
\newcommand{\wrt}{with respect to}

\newcommand{\compose}{\circ}

\DeclareMathOperator{\dist}{dist}

\theoremstyle{plain}
\newtheorem{theorem}{Theorem}[section]
\newtheorem{lemma}[theorem]{Lemma}
\newtheorem{proposition}[theorem]{Proposition}
\theoremstyle{definition}
\newtheorem{definition}[theorem]{Definition}
\newtheorem{example}[theorem]{Example}
\theoremstyle{remark}
\newtheorem{remark}[theorem]{Remark}

\begingroup\makeatletter\ifx\SetFigFont\undefined%
\gdef\SetFigFont#1#2#3#4#5{%
  \reset@font\fontsize{#1}{#2pt}%
  \fontfamily{#3}\fontseries{#4}\fontshape{#5}%
  \selectfont}%
\fi\endgroup%

\begin{document}

\maketitle

\section{Introduction}

Let $(X, \d\mu)$ be a measure space,
$H$ a Hilbert space and $\sigma>0$.

Consider a family of linear operators
$U(t): H \to L_X^2$ defined for each $t \in \R$.
Let $U^*(t): L_X^2 \to H$ be the adjoint of $U(t)$.
We assume that the family $U(t)$
satisfies the \emph{energy estimate}
\begin{equation} \label{eq:energy}
  \norm{U(t) h}_{L_X^2} \les \norm{h}_H, 
  \quad \forall t \in \R, \forall h \in H,
\end{equation}
and the \emph{dispersive inequality}%
\footnote{
  Observe that $L_X^1 \intersection L_X^2$
  is a dense subset of $L_X^1$.}
\begin{equation} \label{eq:dispersion}
  \norm{U(t) U^*(s) f}_{L_X^\infty} \les
  |t-s|^{-\sigma} \norm{f}_{L_X^1},
  \qquad \forall t \ne s,
  \quad \forall f \in L_X^1 \intersection L_X^2.
\end{equation}

The energy estimate allows us to consider the operator
$T: H \to \LtLX{\infty}{2}$ defined as $Th(t) = U(t)h$,
for $t \in \R$ and $h \in H$.
Its formal adjoint is the operator $T^*: \LtLX{1}{2} \to H$
given by the $H$-valued integral
\begin{equation*}
  T^*F = \int U^*(s)F(s) \d s.
\end{equation*}
The composition $TT^*$ is the operator
\begin{equation*}
  TT^*F(t) = \int U(t)U^*(s)F(s) \d s,
\end{equation*}
which can be decomposed as the sum
of its retarded and advanced parts,
\begin{align*}
  (TT^*)_R F(t) &= \int_{s<t} U(t)U^*(s)F(s) \d s, &
  (TT^*)_A F(t) &= \int_{s>t} U(t)U^*(s)F(s) \d s.
\end{align*}
In usual applications, the operator $T$ solves
the initial value problem for
a linear homogeneous differential equation,
while the retarded operator $(TT^*)_R$
solves the corresponding inhomogeneous problem
with zero initial conditions
(Duhamel's principle).

\begin{definition}
  Following~\cite{KeeTao1998},
  we say that the exponent pair $(q,r)$
  is \emph{sharp \mbox{$\sigma$-admissible}} if
  \begin{align*}
    &2 \le q, r \le \infty, &
    &\1q = \sigma \tonde{\12 - \1r}, &
    &(q,r,\sigma) \ne(2,\infty,1).
  \end{align*}
\end{definition}

\begin{definition}
  We introduce another definition
  and say that the pair $(q,r)$ 
  is \emph{$\sigma$-acceptable} if
  \begin{align*}
    &1 \le q,r \le \infty, &
    & \1q < 2\sigma \tonde{\12 - \1r}, &
    & \text{or} \quad (q,r) = (\infty,2).
  \end{align*}
\end{definition}

\begin{figure}
  \centering
  \begin{picture}(0,0)%
    \includegraphics{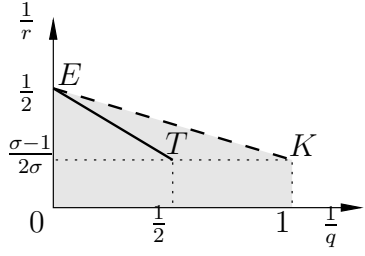}%
  \end{picture}%
  \begin{picture}(2258,1472)(5,-726)
    \put( 151,-661){\makebox(0,0)[lb]{\smash{{\SetFigFont{12}{14.4}{%
              \familydefault}{\mddefault}{\updefault}{$0$}}}}}
    \put(  76, 614){\makebox(0,0)[lb]{\smash{{\SetFigFont{12}{14.4}{%
              \familydefault}{\mddefault}{\updefault}{$\1r$}}}}}
    \put(  76, 164){\makebox(0,0)[lb]{\smash{{\SetFigFont{12}{14.4}{%
              \familydefault}{\mddefault}{\updefault}{$\12$}}}}}
    \put( 324, 258){\makebox(0,0)[lb]{\smash{{\SetFigFont{12}{14.4}{%
              \familydefault}{\mddefault}{\updefault}{$E$}}}}}
    \put( 908,-661){\makebox(0,0)[lb]{\smash{{\SetFigFont{12}{14.4}{%
              \familydefault}{\mddefault}{\updefault}{$\12$}}}}}
    \put(1692,-661){\makebox(0,0)[lb]{\smash{{\SetFigFont{12}{14.4}{%
              \familydefault}{\mddefault}{\updefault}{$1$}}}}}
    \put(1988,-668){\makebox(0,0)[lb]{\smash{{\SetFigFont{12}{14.4}{%
              \familydefault}{\mddefault}{\updefault}{$\1q$}}}}}
    \put(1779,-184){\makebox(0,0)[lb]{\smash{{\SetFigFont{12}{14.4}{%
              \familydefault}{\mddefault}{\updefault}{$K$}}}}}
    \put(   5,-222){\makebox(0,0)[lb]{\smash{{\SetFigFont{12}{14.4}{%
              \familydefault}{\mddefault}{\updefault}{%
              $\frac{\sigma-1}{2\sigma}$}}}}}
    \put(1006,-174){\makebox(0,0)[lb]{\smash{{\SetFigFont{12}{14.4}{%
              \familydefault}{\mddefault}{\updefault}{$T$}}}}}
  \end{picture}%
  \caption{Admissible and acceptable exponents}
  \label{fig:admissible}
\end{figure}

In~\cite{KeeTao1998}, the following theorem was proved:

\begin{theorem}
  If $U(t)$ obeys~\eqref{eq:energy} and~\eqref{eq:dispersion},
  then the estimates
  \begin{align}
    \label{eq:homogest}
    \norm{T h}_{\LtLX{q}{r}} &\les \norm{h}_H, \\
    \label{eq:homogest*}
    \norm{T^* F}_H &\les \norm{F}_{\LtLX{q'}{r'}}, \\
    \label{eq:inhomest}
    \norm{(TT^*)_R F}_{\LtLX{q}{r}} &\les
    \norm{F}_{\LtLX{\qt'}{\rt'}},
  \end{align}
  hold for all sharp $\sigma$-admissible pairs
  $(q,r)$ and $(\qt, \rt)$.
\end{theorem}

As it was already remarked in~\cite{KeeTao1998}, 
we expect the inhomogeneous estimate~\eqref{eq:inhomest}
to have a wider range of admissibility
than the one given by sharp $\sigma$-admissible pairs.
This phenomenon has already been observed
by Harmse~\cite{Har1990} and Oberlin~\cite{Obe1989}
in the context of the wave equation for the case $q=r$.
What they prove is essentially that the inhomogeneous estimate
\begin{equation*}
  \norm{(TT^*)_R F}_{L^p(\R \times X)} \les
  \norm{F}_{L^{\pt'}(\R \times X)},
\end{equation*}
holds when $p, \pt$ satisfy the conditions
\begin{equation*}
  \1p + \1\pt = \frac\sigma{1 + \sigma}, \quad
  \text{and} \quad
  \1p, \1\pt < \frac\sigma{1 + 2\sigma}.
\end{equation*}
Note that the pair $(p,p)$ is sharp $\sigma$-admissible
only for $1/p = \sigma/(2(1+\sigma))$,
while it is $\sigma$-acceptable if and only if $p < \sigma/(2\sigma+1)$.

Also, in the context of Schr\"odinger's equation,
Kato~\cite{Kat1994} proved that
the inhomogeneous estimate~\eqref{eq:inhomest}
holds when the pairs $(q,r)$ and $(\qt, \rt)$
are $\sigma$-acceptable and
satisfy the conditions
\begin{equation*}
  \1q + \1\qt = \sigma \tonde{1 - \1r - \1\rt}, \quad
  \text{and} \quad
  \1r, \1\rt > \frac{\sigma-1}{2\sigma}. 
\end{equation*}

Our goal is to find the largest range
for the pairs $(q,r)$ and $(\qt, \rt)$
which guarantees the validity
of the inhomogeneous estimates~\eqref{eq:inhomest},
and which can be deduced by assuming only 
the energy and dispersive properties,%
~\eqref{eq:energy} and~\eqref{eq:dispersion}.
Our main result is summarized
by the following theorem.

\begin{theorem}[Global inhomogeneous estimates]
  \label{thm:global}
  Let $1 \le q,\qt,r,\rt \le \infty$.
  If $U(t)$ obeys~\eqref{eq:energy} and~\eqref{eq:dispersion},
  then the estimate~\eqref{eq:inhomest} holds
  when the pairs $(q,r)$ and $(\qt,\rt)$ are $\sigma$-acceptable,
  verify the scaling condition
  \begin{equation}
    \label{eq:scaling}
    \1q + \1\qt = \sigma \tonde{1 - \1r - \1\rt}
  \end{equation}
  and satisfy one of the following sets of conditions:
  \begin{itemize}
  \item if $\sigma<1$, there are no further conditions;
  \item if $\sigma=1$, we also require that $r, \rt < \infty$;
  \item if $\sigma>1$, we distinguish two cases,
    \begin{itemize}
    \item non sharp case:
      \begin{align}
        \label{eq:qqt1}
        & \1q + \1\qt < 1, \\
        \label{eq:s1rsr}
        & \frac{\sigma - 1}r \le \frac\sigma{\rt}, \qquad
        \frac{\sigma - 1}\rt \le \frac\sigma{r};
      \end{align}
    \item sharp case:
      \begin{align}
        \label{eq:qqt2}
        & \1q + \1\qt = 1, \\
        \label{eq:s1rsr2}
        & \frac{\sigma - 1}r < \frac\sigma{\rt}, \qquad
        \frac{\sigma - 1}\rt < \frac\sigma{r}, \\
        & \1r \le \1q, \qquad \1\rt \le \1\qt.
      \end{align}
    \end{itemize}
  \end{itemize}
\end{theorem}

\begin{remark}
  Conditions~\eqref{eq:qqt1} and~\eqref{eq:s1rsr}
  which appear in the non sharp case for $\sigma>1$
  are always trivially satisfied if $\sigma<1$
  or if $\sigma=1$ and $r,\rt <\infty$.
\end{remark}

\begin{remark}
  Condition~\eqref{eq:scaling} together with $1/q+1/\qt \le 1$
  have the following interpretation:
  if $(1/Q, 1/R)$ is the midpoint between
  the points $(1/q, 1/r)$ and $(1/\qt, 1/\rt)$,
  then $(Q,R)$ is a sharp $\sigma$-admissible pair.
\end{remark}

\begin{remark}
  Formally, it is easy to verify that $TT^*$
  coincides with its dual $(TT^*)^*$,
  while ${(TT^*)_R}^* = (TT^*)_A$.
  Moreover, $(TT^*)_A$ becomes $(TT^*)_R$
  if we invert the direction of time.
  These duality relations explain why all conditions
  must be invariant under the symmetry
  $(q,r) \leftrightarrow (\qt, \rt)$.
\end{remark}

\begin{remark}
  In the case $q=r$ and $\qt=\rt$,
  theorem~\ref{thm:global}
  reduces to the results of Harmse and Oberlin
  which can be shown to be optimal.
\end{remark}

\begin{figure}
  \centering
  \begin{picture}(0,0)%
    \includegraphics{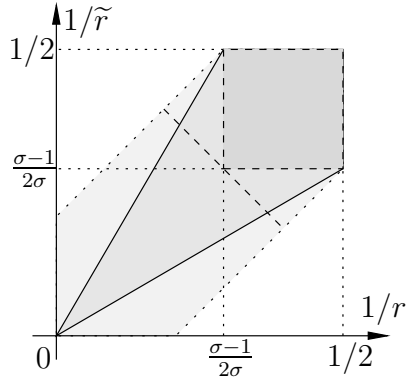}%
  \end{picture}%
  \begin{picture}(2596,2365)(0,-1514)
    \put(1999,-1452){\makebox(0,0)[lb]{\smash{{\SetFigFont{12}{14.4}{%
              \familydefault}{\mddefault}{\updefault}{$1/2$}}}}}
    \put( 172,-1465){\makebox(0,0)[lb]{\smash{{\SetFigFont{12}{14.4}{%
              \familydefault}{\mddefault}{\updefault}{$0$}}}}}
    \put(   0,  475){\makebox(0,0)[lb]{\smash{{\SetFigFont{12}{14.4}{%
              \familydefault}{\mddefault}{\updefault}{$1/2$}}}}}%
    \put(2212,-1145){\makebox(0,0)[lb]{\smash{{\SetFigFont{12}{14.4}{%
              \familydefault}{\mddefault}{\updefault}{$1/r$}}}}}
    \put(1245,-1445){\makebox(0,0)[lb]{\smash{{\SetFigFont{12}{14.4}{%
              \familydefault}{\mddefault}{\updefault}{%
              $\frac{\sigma-1}{2\sigma}$}}}}}
    \put( 345,  662){\makebox(0,0)[lb]{\smash{{\SetFigFont{12}{14.4}{%
              \familydefault}{\mddefault}{\updefault}{$1/\rt$}}}}}
    \put(  12, -258){\makebox(0,0)[lb]{\smash{{\SetFigFont{12}{14.4}{%
              \familydefault}{\mddefault}{\updefault}{%
              $\frac{\sigma-1}{2\sigma}$}}}}}
  \end{picture}%
  \caption{Admissible range for the exponents $r$ and $\rt$.}
  \label{fig:rtr}
\end{figure}

\begin{remark}
  When $\sigma>1$,
  theorem~\ref{thm:global} improves on Kato's result~\cite{Kat1994}.
  Kato's theorem required $r$ and $\rt$ to be less than $2\sigma/(\sigma-1)$.
  We replace that restriction with a condition
  which can be read as
  \begin{equation*}
    \frac{\sigma-1}{\sigma} \le \frac{r}{\rt}
    \le \frac{\sigma}{\sigma-1}.
  \end{equation*}
\end{remark}

Our proof of theorem~\ref{thm:global}
makes use of the techniques of Keel and Tao~\cite{KeeTao1998}
and is based on the following localized version 
of the inhomogeneous estimates.

\begin{theorem}[Local inhomogeneous estimates]
  \label{thm:local}
  Assume $U(t)$ obeys~\eqref{eq:energy} and~\eqref{eq:dispersion},
  and let $I$ and $J$ be two time intervals 
  of unit length $\abs{I} = \abs{J} = 1$
  separated by a distance of scale $1$, $\dist(I,J) \sss 1$.
  Then, the estimate
  \begin{equation}
    \label{eq:localest}
    \norm{TT^* F}_{\LtLX[J]{q}{r}}
    \les \norm{F}_{\LtLX[I]{\qt'}{\rt'}}, \qquad
    \forall F \in \LtLX[I]{\qt'}{\rt'},
  \end{equation}
  holds for all pairs $(q,r)$ and $(\qt, \rt)$
  which verify the following conditions:
  \begin{gather}
    \label{eq:clie1}
    1 \le q, \qt \le \infty, \qquad
    2 \le r, \rt \le \infty, \\
    \label{eq:clie2}
    \frac{\sigma - 1}r \le \frac\sigma{\rt}, \qquad
    \frac{\sigma - 1}\rt \le \frac\sigma{r}, \\
    \label{eq:clie3}
    \1q \ge \sigma \tonde{\1\rt - \1r}, \qquad
    \1\qt \ge \sigma \tonde{\1r - \1\rt},
  \end{gather}
  and if $\sigma = 1$, we must also require $r, \rt < \infty$.
\end{theorem}

\begin{remark}
  For the local estimates of theorem~\ref{thm:local}
  we do not require the pairs $(q,r)$ and $(\qt,\rt)$
  to be $\sigma$-acceptable.
\end{remark}

\section{Proof of the local estimates}
\label{sec:proof-local-estim}

\begin{proof}[Proof of theorem~\ref{thm:local}.]
  Let $\Elocal$ be the set of points
  $\tonde{1/q, 1/r; 1/\qt, 1/\rt}$ in $[0,1]^4$
  corresponding to the pairs $(q,r)$, $(\qt, \rt)$
  for which the estimate~\eqref{eq:localest} is valid.
  
  We start by observing that the
  dispersive estimate~\eqref{eq:dispersion}
  immediately yields the case $q=r=\qt=\rt=\infty$,
  \begin{multline}
    \label{eq:0000}
    \norm{TT^* F}_{\LtLX[J]{\infty}{\infty}} \le
    \int_I\norm{U(t)U^*(s)F(s)}_{\LtLX[J]{\infty}{\infty}}\d s
    \les \\ \les \int_I \norm{F(s)}_{L_X^1} \d s
    = \norm{F}_{\LtLX[I]{1}{1}}.  
  \end{multline}
  Hence, $(0,0;0,0) \in \Elocal$.
  
  On the other hand, if we exploit the factorization $TT^*$,
  we can apply the homogeneous Strichartz estimates%
  ~\eqref{eq:homogest} and~\eqref{eq:homogest*},
  \begin{equation}
    \label{eq:qrqtrt}
    \norm{TT^* F}_{\LtLX[J]{q}{r}} \les
    \norm{T^* F}_H \les 
    \norm{F}_{\LtLX[I]{\qt'}{\rt'}},
  \end{equation}
  and obtain that
  $\tonde{1/q, 1/r; 1/\qt, 1/\rt} \in \Elocal$
  whenever $(q,r)$ and $(\qt, \rt)$
  are sharp $\sigma$-admissible pairs.
  
  By standard $L^p$ interpolation%
  \footnote{See~\cite{BerLof1976}.}
  between~\eqref{eq:0000} and~\eqref{eq:qrqtrt},
  we obtain that $\Elocal$ contains
  the convex hull of the set 
  \begin{equation*}
    \graffe{\tonde{0, 0; 0, 0}} \union
    \graffe{\tonde{\1q, \1r; \1\qt, \1\rt}:
      \text{$(q,r)$ and $(\qt, \rt)$
        are sharp $\sigma$-admissible pairs}}.
  \end{equation*}
  
  Since we have restricted $F$ and $TT^*F$ to unit time intervals,
  it follows from H\"older's inequality
  that when $q \ge Q$, $\qt \ge \Qt$ and
  $\tonde{\1q,\1r; \1\qt,\1\rt} \in \Elocal$ then
  $\tonde{\1Q,\1r; \1\Qt,\1\rt} \in \Elocal$.
  If we apply this property
  to the points of the above convex hull
  we obtain that $\Elocal$ contains a set $\mE_*$
  exactly described by the conditions
  appearing in theorem~\ref{thm:local}.
  More details of this computation
  are given in appendix~\ref{sec:details-about-set}.
\end{proof}  

\begin{remark}
  As observed in~\cite{KeeTao1998},
  there exists a natural scaling
  associated to the family $U(t)$.
  More precisely, let $\lambda>0$, then
  the hypotheses~\eqref{eq:energy} and~\eqref{eq:dispersion}
  are invariant under the rescaling
  \begin{align*}
    & \d\mu \leftarrow \lambda^\sigma \d\mu, &
    & \ip{h_1,h_2}_H \leftarrow
    \lambda^\sigma \ip{h_1,h_2}_H,&
    & U(t) \leftarrow U(t/\lambda).
  \end{align*}
  If we apply theorem~\ref{thm:local}
  to the rescaled operators
  and espress the result in terms of the original operators
  we obtain the following generalization
  of the local estimates.
\end{remark}

\begin{proposition} \label{pro:local-scaled}
  Let $I$ and $J$ be two time intervals 
  of length $\lambda$,
  $\abs{I} = \abs{J} = \lambda$,
  separated by a distance of scale $\lambda$,
  $\dist(I,J) \sss \lambda$.
  We have the estimate
  \begin{equation}
    \label{eq:scaledest}
    \norm{TT^* F}_{\LtLX[J]{q}{r}}
    \les \lambda^{\beta(q,r;\qt,\rt)}
    \norm{F}_{\LtLX[I]{\qt'}{\rt'}},
  \end{equation}
  with
  \begin{equation*}
    \beta(q,r;\qt,\rt) =
    \1q + \1\qt - \sigma \tonde{1 - \1r - \1\rt},
  \end{equation*}
  whenever the pairs $(q,r)$ and $(\qt, \rt)$
  satisfy the conditions appearing in theorem~\ref{thm:local}.
\end{proposition}

\section{Dyadic decompositions
  of sets, functions and operators}

By duality, the linear estimate~\eqref{eq:inhomest}
is equivalent to the bilinear estimate
\begin{equation*}
  \abs{B(F,G)} \le
  \norm{F}_{\LtLX{\qt'}{\rt'}} \norm{G}_{\LtLX{q'}{r'}},
\end{equation*}
where $B$ is the scalar bilinear operator
\begin{equation} \label{eq:BFG}
  B(F,G) = \iint_{s<t}
  \ip{U(s)^* F(s), U(t)^* G(t)} \d s \d t.
\end{equation}

We want to decompose $B$ into a sum of localized operators
to which we can apply proposition~\ref{pro:local-scaled}.
In order to do so, we make use of
Whitney's decompositions of open sets
applied to the domain of the integration in~\eqref{eq:BFG}.

We say that $\lambda$ is dyadic number if $\lambda = 2^k$
for some integer $k$.
The set of all dyadic numbers, $2^\Z$,
is a multiplicative abelian group.
In the following $\lambda$, $\mu$ and $\nu$
will always denote dyadic numbers.
In particular, if $\alpha>0$ then
\begin{equation*}
  \sum_{\lambda: \lambda<\mu} \lambda^\alpha = \1{2^\alpha-1} \cdot \mu^\alpha
\end{equation*}
is the sum of a convergent geometric series.

Recall that a \emph{dyadic square} in $\R^2$
is a square whose sidelength
is a dyadic number $\lambda \in 2^\Z$
and such that all the coordinates of its vertices
are integer multiples of $\lambda$.

\begin{lemma}[Dyadic Whitney decomposition (see~\cite{Ste1970})]
  \label{lem:whitney}
  Let $\Omega$ be a proper open subset of $\R^2$.
  There exists a partition of $\Omega$
  into a family $\mQ$ of essentially disjoint%
  \footnote{By essentially disjoint we mean that
    the interiors of the squares are disjoint.}
  dyadic squares
  with the property that the distance of $Q$
  from the boundary of $\Omega$
  is approximately proportional to the diameter of $Q$.
\end{lemma}

\begin{figure}
  \centering
  \includegraphics{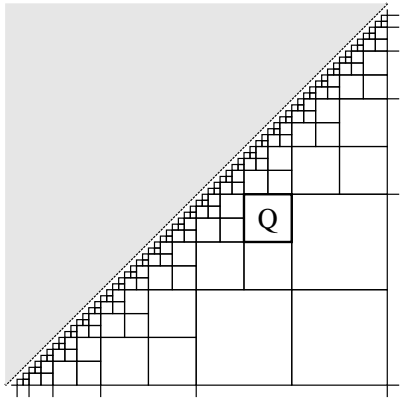}
  \caption{Whitney's decomposition of the region $s<t$.}
  \label{fig:whitney}
\end{figure}

Let $\mQ$ be the Whitney decomposition
for the domain $\Omega = \graffe{(s,t): s<t}$
given by lemma~\ref{lem:whitney}.
For each dyadic number $\lambda$,
let $\mQ_\lambda$ be the family of squares in $\mQ$
whose sidelength is $\lambda$.
Each square $Q = I \times J \in \mQ_\lambda$,
has the property that
\begin{equation}
  \label{eq:IJQ}
  \lambda = \abs{I} = \abs{J} \sss \dist(Q, \de\Omega) \sss \dist(I,J).
\end{equation}

Since $\Omega = \union_\lambda \union_{Q \in \mQ_\lambda} Q$
and the squares $Q$ are essentially disjoint,
we can write the decomposition
\begin{equation}
  \label{eq:BBQ}
  B = \sum_\lambda \sum_{Q \in \mQ_\lambda} B_Q,
\end{equation}
where, for each square $Q = I \times J$, we set
\begin{equation*}
  B_Q(F,G) = B(\chi_I F, \chi_J G) = 
  \iint_{\substack{s \in I \\ t \in J}}
  \ip{U(s)^* F(s), U(t)^* G(t)} \d s \d t,
\end{equation*}
with $\chi_I$, $\chi_J$ being the characteristic functions
of the intervals $I$ and $J$.
The local estimate~\eqref{eq:scaledest} of proposition~\ref{pro:local-scaled}
is equivalent to a bilinear estimate
for the localized operator $B_Q$, $Q=I \times J$, namely
\begin{equation}
  \label{eq:BQFGest}
  \abs{B_Q(F,G)} \les \lambda^{\beta(q,r;\qt,\rt)} 
  \norm{F}_{\LtLX[I]{\qt'}{\rt'}} \norm{G}_{\LtLX[J]{q'}{r'}}.
\end{equation}

\begin{lemma} \label{lem:sumQIJ}
  Suppose that $1/p + 1/\pt \ge 1$, then we have
  \begin{equation*}
    \sum_{\substack{Q \in \mQ_\lambda \\ Q = I \times J}}
    \norm{f}_{L^\pt(I)} \norm{g}_{L^p(J)}  \le
    \norm{f}_{L^\pt(\R)} \norm{g}_{L^p(\R)},
  \end{equation*}
  for $f \in L^\pt(\R)$, $g \in L^p(\R)$ and any dyadic number $\lambda$.
\end{lemma}

\begin{proof}
  It follows immediately from the inequality
  \begin{equation*}
    \sum_n \abs{A_n B_n} \le
    \TOnde{\sum_n {\abs{A_n}^\pt}}^{1/\pt}
    \TOnde{\sum_n {\abs{B_n}^p}}^{1/p},
  \end{equation*}
  which is valid only if $1/p + 1/\pt \ge 1$,
  and the fact that for each dyadic interval $I$
  there are at most a fixed finite number of intervals $J$
  which satisfy~\eqref{eq:IJQ} and they are all contained
  in a neighborhood of $I$ of size $O(\lambda)$.
\end{proof}

It follows from lemma~\ref{lem:sumQIJ}
that if~\eqref{eq:BQFGest} holds with $1/q + 1/\qt \le 1$
then we have
\begin{equation}
  \label{eq:sumQBFG}
  \sum_{Q \in \mQ_\lambda} \abs{B_Q(F,G)} \les \lambda^{\beta(q,r;\qt,\rt)} 
  \norm{F}_{\LtLX{\qt'}{\rt'}} \norm{G}_{\LtLX{q'}{r'}}.
\end{equation}

As in Keel and Tao~\cite{KeeTao1998},
we are going to decompose the functions $F$ and $G$
into dyadic atoms
and then play with interpolation
on the bilinear version of our operators.
There is a small difference
with respect to the approach of Keel and Tao:
in the non sharp case of theorem~\ref{thm:global},
instead of decomposing functions, for each fixed time,
into dyadic pieces with respect to the $L_X^r$ norm,
we are going to decompose our functions
into $L_X^r$-valued dyadic pieces
with respect to the $L_t^q$ norm.
This will allow us to recover some extreme cases,
namely the cases of equality in condition~\eqref{eq:s1rsr}.
For the sharp case of theorem~\ref{thm:global},
we will perform the same dyadic decomposition as in Keel and Tao.

Let $1 \le p \le \infty$.
Let $\mX$ be a measurable space and $\mB$ be a Banach space.
A \emph{$p$-atom} in $L^p(\mX;\mB)$ of size $\lambda$
is a measurable function $\varphi:\mX \to \mB$
such that
\begin{itemize}
\item
  $\xi \mapsto \varphi(\xi)$ is supported
  on a set of measure less than $\lambda$;
\item
  $\norm{\varphi}_{L^\infty(\mX; \mB)} \les \lambda^{1/p}$.
\end{itemize}
It follows that 
we have $\norm{\varphi}_{L^p(\mX; \mB)} \les 1$
for any $p$-atom $\varphi$.
More generally, for any $p$-atom of size $\lambda$
and any exponent $q \in [1,\infty]$ we have
\begin{equation}
  \label{eq:atomqp}
  \norm{\varphi}_{L^q(\mX;\mB)} \les
  \lambda^{\frac1q - \frac1p}.
\end{equation}

Any $L^p$ function can be decomposed into a dyadic sum
of $p$-atoms.

\begin{lemma} \label{lem:atomic}
  Any $\mB$-valued function $F \in L^p(\mX;\mB)$
  can be decomposed as
  \begin{equation*}
    F(t) =
    \sum_{\lambda \in 2^\Z} a_\lambda \varphi_\lambda(t),
  \end{equation*}
  where
  \begin{itemize}
  \item each $\varphi_\lambda(t)$
    is a $p$-atom in $L^p(\mX; \mB)$ of size $\lambda$;
  \item the atoms $\varphi_\lambda$ have disjoint supports;
  \item $a_\lambda$ are non-negative constants such that
    \begin{equation*}
      \norm{f}_{L^p(\mX; \mB)} \sss
      \norm{a_\lambda}_{\ell^p\tonde{2^\Z}}.
    \end{equation*}
  \end{itemize}
\end{lemma}

The proof of the lemma is the same as in the scalar case
(see lemma $5.1$ in~\cite{KeeTao1998}).

\section{Proof of the global estimates: non sharp case.}

We assume now that we are in the non sharp case with $1/q + 1/\qt < 1$
and we want to prove
the global inhomogeneous Strichartz estimates~\eqref{eq:inhomest}
under the conditions stated in theorem~\ref{thm:global}.

We can apply lemma~\ref{lem:atomic} to the functions
$F \in \LtLX{\qt'}{\rt'}$ and $G \in \LtLX{q'}{r'}$,
and obtain the decompositions
\begin{align*}
  & F(t) = \sum_\mu a_{\mu} \varphi_{\mu}(t), &
  & G(t) = \sum_\nu b_{\nu} \psi_{\nu}(t),
\end{align*}
where $\varphi_\mu$ is a $\qt'$-atom
in $\LtLX{\qt'}{\rt'}$ of size $\mu$,
$\psi_\nu$ is a $q'$-atom
in $\LtLX{q'}{r'}$ of size $\nu$
and
\begin{align}
  \label{eq:FaGb}
  & \norm{F}_{\LtLX{\qt'}{\rt'}} \sss
  \norm{a_\mu}_{\ell^{\qt'}}, &
  & \norm{G}_{\LtLX{q'}{r'}} \sss
  \norm{b_\nu}_{\ell^{q'}}.
\end{align}
We plug these decompositions into our previous decomposition~\eqref{eq:BBQ}
of the operator $B$ and obtain
\begin{equation}
  \label{eq:BFGsum}
  B(F,G) = \sum_{\lambda, \mu, \nu} a_\mu b_\nu
  \sum_{Q \in \mQ_\lambda} B_Q(\varphi_\mu, \psi_\nu).
\end{equation}

If we apply the local estimates for the terms $B_Q(\varphi_\mu, \psi_\nu)$
directly to this sum we will obtain a divergent sum.
Fortunately, as was well illustrated in~\cite{KeeTao1998},
we can gain some summability by slightly perturbing
the exponents $\qt$ and $q$.

\begin{remark}
  To simplify the notation it is convenient
  to introduce the function $[\cdot]: \R^+ \to \R^+$ defined by
  \begin{equation*}
    [\lambda] = \max\GRaffe{\lambda, \1\lambda},
  \end{equation*}
  which plays the role, in the multiplicative group $\R^+$,
  played by the absolute value in the additive group $\R$.
  In particular $[2^k] = 2^{|k|}$.
\end{remark}

\begin{lemma} \label{lem:epsiqq}
  Suppose $(q_0,r)$ and $(\qt_0,\rt)$ are such that
  $1/q_0 + 1/\qt_0 < 1$ and that
  the estimate~\eqref{eq:sumQBFG} holds
  with exponents $(q,r)$ and $(\qt,\rt)$
  for all $(1/q, 1/\qt)$ in a full neighborhood of $(1/q_0, 1/\qt_0)$.
  Then, there exists $\epsi>0$ such that,
  for all dyadic numbers $\lambda, \mu, \nu$, we have
  \begin{equation*}
    \sum_{Q \in \mQ_\lambda} \abs{B_Q(\varphi_\mu, \psi_\nu)} \les
    \lambda^{\beta(q_0,r;\qt_0,\rt)}
    \quadre{\frac\mu\lambda}^{-\epsi} \quadre{\frac\nu\lambda}^{-\epsi},
  \end{equation*}
  whenever $\varphi_\mu$ is a $\qt_0'$-atom of size $\mu$
  in $\LtLX{\qt_0'}{\rt'}$,
  and $\psi_\nu$ is a $q_0'$-atom in $\LtLX{q_0'}{r'}$ of size $\nu$.
\end{lemma}

\begin{proof}
  When $(1/q, 1/\qt)$ is close enough to $(1/q_0, 1/\qt_0)$
  we still have $1/q + 1/\qt <1$.
  Hence, we can combine the estimate~\eqref{eq:sumQBFG}
  with the property~\eqref{eq:atomqp} of dyadic atoms
  and we obtain
  \begin{equation*}
    \sum_{Q \in \mQ_\lambda} \abs{B_Q(\varphi_\mu, \psi_\nu)} \les
    \lambda^{\beta(q,r;\qt,\rt)} \mu^{\1{\qt_0}-\1\qt} \nu^{\1{q_0}-\1q} =
    \lambda^{\beta(q_0,r;\qt_0,\rt)}
    \tonde{\frac\mu\lambda}^{\1{\qt_0}-\1\qt}
    \tonde{\frac\nu\lambda}^{\1{q_0}-\1q}.
  \end{equation*}
  For given $\lambda, \mu, \nu$, we can choose $\qt$ and $q$
  in a neighborhood of $\qt_0$ and $q_0$ so that 
  \begin{align*}
    & \1{\qt_0}-\1\qt = 
    \begin{sistema}
      & +\epsi, & \text{if $\mu \le \lambda$,} \\
      & -\epsi, & \text{if $\mu > \lambda$;}
    \end{sistema} &
    & \1{q_0}-\1q = 
    \begin{sistema}
      & +\epsi, & \text{if $\nu \le \lambda$,} \\
      & -\epsi, & \text{if $\nu > \lambda$,}
    \end{sistema}
  \end{align*}
  where $\epsi$ is a small positive number
  (independent of $\lambda, \mu, \nu$).
  Doing this way we have
  \begin{align*}
    & \tonde{\frac\mu\lambda}^{\1{\qt_0}-\1\qt} =
    \quadre{\frac\mu\lambda}^{-\epsi}, &
    & \tonde{\frac\nu\lambda}^{\1{q_0}-\1q} =
    \quadre{\frac\nu\lambda}^{-\epsi},     
  \end{align*}
  and the lemma is proved.
\end{proof}

The inequality~\eqref{eq:sumQBFG} holds
when $(q,r)$ and $(\qt,\rt)$
are in the range of validity of the local estimate~\eqref{eq:localest},
described by conditions~\eqref{eq:clie1},~\eqref{eq:clie2},~\eqref{eq:clie3}.
In order to apply lemma~\ref{lem:epsiqq}
to the sum in~\eqref{eq:BFGsum}
we require that $\1q + \1\qt < 1$
and that we have strict inequalities in~\eqref{eq:clie3},
so that~\eqref{eq:sumQBFG} remains valid
under small perturbations of $q$ and $\qt$.
We obtain 
\begin{equation*}
  \abs{B(F,G)}\le \sum_{\mu, \nu} a_\mu b_\nu
  \sum_\lambda \lambda^{\beta(q,r;\qt,\rt)}
    \quadre{\frac\mu\lambda}^{-\epsi} \quadre{\frac\nu\lambda}^{-\epsi}.
\end{equation*}
The sum over $\lambda$ diverges unless $\beta(q,r;\qt,\rt) = 0$,
in which case we have
\begin{equation*}
  \sum_\lambda
  \quadre{\frac\mu\lambda}^{-\epsi} \quadre{\frac\nu\lambda}^{-\epsi} \les
  \tonde{1 + \log\quadre{\frac\mu\nu}} \quadre{\frac\mu\nu}^{-\epsi} = 
  c_{\mu/\nu},
\end{equation*}
where the sequence $c_\lambda = (1 + \log[\lambda]) [\lambda]^{-\epsi}$
is summable,
\begin{equation*}
  \norm{c_\lambda}_{\ell^1} =
  \sum_\lambda (1 + \log[\lambda]) [\lambda]^{-\epsi} = 
  \sum_{n \in \Z} (1+|n|) 2^{-\epsi|n|} < \infty.
\end{equation*}
Hence,
\begin{equation}
  \label{eq:BFGabc}
  \abs{B(F,G)} \les \sum_{\mu, \nu} a_\mu b_\nu c_{\mu/\nu},
\end{equation}
where the right hand side
is a convolution sum written in multiplicative index notation.

\begin{lemma}[Young's inequality for convolution of sequences]
  \label{lem:Young}
  Let $A_n$, $B_n$, $C_n$ be sequences of non negative numbers.
  If
  \begin{equation*}
    \1p + \1q + \1r \ge 2,
  \end{equation*}
  then $\sum_{n,k} A_n B_k C_{n-k} \le 
    \norm{A}_{\ell^p} \norm{B}_{\ell^q} \norm{C}_{\ell^r}$.
\end{lemma}

Since we have
\begin{equation*}
  \1{\qt'} + \1{q'} + 1 = 3 - \1\qt - \1q > 2,
\end{equation*}
we can apply Young's inequality to~\eqref{eq:BFGabc}
and use~\eqref{eq:FaGb} to finally obtain
\begin{equation*}
  \abs{B(F,G)} \les 
  \norm{a}_{\ell^{\qt'}} \norm{b}_{\ell^{q'}} \norm{c}_{\ell^1} \les
  \norm{F}_{\LtLX{\qt'}{\rt'}} \norm{G}_{\LtLX{q'}{r'}}.
\end{equation*}

We summarize the conditions we have imposed so far
on the parameters $q,r,\qt,\rt$:
\begin{itemize}
\item the non sharp case condition $1/q + 1/\qt < 1$;
\item the scaling invariant condition $\beta(q,r;\qt,\rt) = 0$,
  which is equivalent to~\eqref{eq:scaling};
\item conditions on $r$ and $\rt$ for the validity of the local estimates,
  \begin{align*}
    & \frac{\sigma - 1}r \le \frac\sigma{\rt}, &
    & \frac{\sigma - 1}\rt \le \frac\sigma{r};
  \end{align*}
\item conditions on $q$ and $\qt$ for the validity of the local estimates
  with strict inequality,
  \begin{align*}
    & \1q > \sigma \tonde{\1\rt - \1r}, &
    & \1\qt > \sigma \tonde{\1r - \1\rt},
  \end{align*}
  which, in the scaling invariant case $\beta(q,r;\qt,\rt) = 0$,
  become equivalent to say that
  $(q,r)$ and $(\qt,\rt)$ are $\sigma$-acceptable pairs.
\end{itemize}

\section{Proof of the global estimates: sharp case.}

In the sharp case of theorem~\ref{thm:global},
if we slightly perturb the values of $q$ and $\qt$,
we may violate the condition $1/q + 1/\qt \le 1$
which is necessary for lemma~\ref{lem:sumQIJ},
and we could not repeat the trick used in lemma~\ref{lem:epsiqq}
to gain summability in~\eqref{eq:BFGsum}.
However, we can still play with the exponents $r$ and $\rt$,
if we perform the atomic decomposition on the function $F(t)$ and $G(t)$
for each fixed $t$.

We assume now that $1/q + 1/\qt = 1$.

For each $t$, we can apply lemma~\ref{lem:atomic} to the functions
$F(t) \in L_X^{\rt'}$ and $G(t) \in L_X^{r'}$,
and obtain the decompositions
\begin{align*}
  & F(t) = \sum_\mu a_{\mu}(t) \varphi_{\mu}(t), &
  & G(t) = \sum_\nu b_{\nu}(t) \psi_{\nu}(t),
\end{align*}
where $\varphi_\mu(t)$ is a $\rt'$-atom
in $L_X^{\rt'}$ of size $\mu$,
$\psi_\nu(t)$ is a $r'$-atom
in $L_X^{r'}$ of size $\nu$
and
\begin{align}
  \label{eq:FtatGtbt}
  & \norm{F(t)}_{L_X{\rt'}} \sss
  \norm{a_\mu(t)}_{\ell^{\rt'}}, &
  & \norm{G(t)}_{L_X{r'}} \sss
  \norm{b_\nu(t)}_{\ell^{r'}}.
\end{align}
We plug these decompositions into our previous decomposition~\eqref{eq:BBQ}
of the operator $B$ and obtain
\begin{equation}
  \label{eq:BFGsum2}
  B(F,G) = \sum_{\lambda, \mu, \nu} 
  \sum_{Q \in \mQ_\lambda} B_Q(a_\mu \varphi_\mu, b_\nu \psi_\nu).
\end{equation}

\begin{lemma} \label{lem:epsirr}
  Suppose $(q,r_0)$ and $(\qt,\rt_0)$ are such that
  the estimate~\eqref{eq:sumQBFG} holds
  with exponents $(q,r)$ and $(\qt,\rt)$
  for all $(1/r, 1/\rt)$ in a full neighborhood of $(1/r_0, 1/\rt_0)$.
  Then, there exists $\epsi>0$ such that,
  for all dyadic numbers $\lambda, \mu, \nu$
  and dyadic square $Q = I \times J \in \mQ_\lambda$, we have
  \begin{equation*}
    \abs{B_Q(a \varphi_\mu, b \psi_\nu)} \les
    \lambda^{\beta(q,r_0;\qt,\rt_0)}
    \quadre{\frac\mu{\lambda^\sigma}}^{-\epsi}
    \quadre{\frac\nu{\lambda^\sigma}}^{-\epsi}
    \norm{a}_{L^{\qt'}(I)} \norm{b}_{L^{q'}(J)},
  \end{equation*}
  whenever $a \in L^{\qt'}(I;\R)$, $b \in L^{q'}(J;\R)$,
  and for each $t$, the function
  $\varphi_\mu(t)$ is a $\rt_0'$-atom of size $\mu$ in $L_X{\rt_0'}$,
  and the function $\psi_\nu(t)$ is a $r_0'$-atom in $L_X{r_0'}$ of size $\nu$.
\end{lemma}

\begin{proof}
  We combine the local estimate~\eqref{eq:BQFGest}
  with property~\eqref{eq:atomqp} of dyadic atoms
  and we obtain
  \begin{multline*}
    \abs{B_Q(\varphi_\mu, \psi_\nu)} \les
    \lambda^{\beta(q,r;\qt,\rt)} 
    \norm{a}_{L^{\qt'}(I)}
    \norm{\varphi_\mu}_{\LtLX[I]{\infty}{\rt'}}
    \norm{b}_{L^{q'}(J)}
    \norm{\psi_\nu}_{\LtLX[J]{q'}{r'}} \les \\
    \les \lambda^{\beta(q,r;\qt,\rt)} \mu^{\1{\rt_0}-\1\rt} \nu^{\1{r_0}-\1r} =
    \lambda^{\beta(q,r_0;\qt,\rt_0)}
    \tonde{\frac\mu{\lambda^\sigma}}^{\1{\rt_0}-\1\rt}
    \tonde{\frac\nu{\lambda^\sigma}}^{\1{r_0}-\1r}.
  \end{multline*}
  Similarly to what we did in the proof of lemma~\ref{lem:epsiqq},
  for given $\lambda, \mu, \nu$, we can choose $\rt$ and $r$
  in a neighborhood of $\rt_0$ and $r_0$ so that 
  \begin{align*}
    & \tonde{\frac\mu{\lambda^\sigma}}^{\1{\rt_0}-\1\rt} =
    \quadre{\frac\mu{\lambda^\sigma}}^{-\epsi}, &
    & \tonde{\frac\nu{\lambda^\sigma}}^{\1{r_0}-\1r} =
    \quadre{\frac\nu{\lambda^\sigma}}^{-\epsi},     
  \end{align*}
  and the lemma is proved.
\end{proof}

The inequality~\eqref{eq:BQFGest} holds
when $(q,r)$ and $(\qt,\rt)$
are in the range of validity of the local estimate~\eqref{eq:localest},
described by conditions~\eqref{eq:clie1},~\eqref{eq:clie2},~\eqref{eq:clie3}.
In order to apply lemma~\ref{lem:epsiqq}
to the sum in~\eqref{eq:BFGsum2},
we require strict inequalities
in~\eqref{eq:clie2} and~\eqref{eq:clie3},
so that~\eqref{eq:BQFGest} remains valid
under small perturbations of $r$ and $\rt$.
We obtain
\begin{equation*}
  \abs{B(F,G)}\le \sum_{\lambda, \mu, \nu}
  \lambda^{\beta(q,r;\qt,\rt)}
  \quadre{\frac\mu{\lambda^\sigma}}^{-\epsi}
  \quadre{\frac\nu{\lambda^\sigma}}^{-\epsi}
  \sum_{I \times J \in \mQ_\lambda} 
  \norm{a_\mu}_{L^{\qt'}(I)} \norm{b_\nu}_{L^{q'}(J)},
\end{equation*}
Since $1/q + 1/\qt = 1$ we can still apply lemma~\ref{lem:sumQIJ}
and the sum reduces to
\begin{equation*}
  \abs{B(F,G)}\le \sum_{\mu, \nu}
  \tonde{\sum_\lambda \lambda^{\beta(q,r;\qt,\rt)}
    \quadre{\frac\mu{\lambda^\sigma}}^{-\epsi}
    \quadre{\frac\nu{\lambda^\sigma}}^{-\epsi}}
  \norm{a_\mu}_{L^{\qt'}(\R)} \norm{b_\nu}_{L^{q'}(\R)}.
\end{equation*}
As in the previous section,
the sum over $\lambda$ diverges unless $\beta(q,r;\qt,\rt) = 0$,
in which case we have
\begin{equation*}
  \sum_\lambda
  \quadre{\frac\mu{\lambda^\sigma}}^{-\epsi}
  \quadre{\frac\nu{\lambda^\sigma}}^{-\epsi} \les
  \tonde{1 + \log\quadre{\frac\mu\nu}} \quadre{\frac\mu\nu}^{-\epsi} = 
  c_{\mu/\nu}.
\end{equation*}
Hence,
\begin{equation}
  \label{eq:BFGABc}
  \abs{B(F,G)}\le \sum_{\mu, \nu}
  \norm{a_\mu}_{L^{\qt'}(\R)} \norm{b_\nu}_{L^{q'}(\R)} c_{\mu/\nu}.
\end{equation}
Since we have
\begin{equation*}
  \1{\qt'} + \1{q'} + 1 = 3 - \1\qt - \1q = 2,
\end{equation*}
we can apply lemma~\ref{lem:Young} to~\eqref{eq:BFGABc}
and obtain
\begin{multline*}
  \abs{B(F,G)} \les
  \TOnde{\sum_\mu \norm{a_\mu(t)}_{L^{\qt'}}^{\qt'}}^{1/{\qt'}}
  \TOnde{\sum_\nu \norm{b_\nu(t)}_{L^{q'}}^{q'}}^{1/{q'}} = \\
  = \NOrm{\TOnde{\sum_\mu a_\mu(t)^{\qt'}}^{1/{\qt'}}}_{L^{\qt'}}
  \NOrm{\TOnde{\sum_\nu b_\nu(t)^{q'}}^{1/{q'}}}_{L^{q'}}.
\end{multline*}
To finish the proof, we observe that we have
\begin{align*}
  & \TOnde{\sum_\mu a_\mu(t)^{\qt'}}^{1/{\qt'}} \le
  \TOnde{\sum_\mu a_\mu(t)^{\rt'}}^{1/{\rt'}} \sss \norm{F(t)}_{L_X{\rt'}}, \\
  & \TOnde{\sum_\nu b_\nu(t)^{q'}}^{1/{q'}} \le
  \TOnde{\sum_\nu b_\nu(t)^{r'}}^{1/{r'}} \sss \norm{G(t)}_{L_X{r'}},
\end{align*}
if we require%
\footnote{
  This seems to be a technical condition 
  which made us prefer to proceed with the different dyadic decomposition
  in the non sharp case,
  but which we are not able to avoid in the sharp case.}
that $\qt \le \rt$ and $q \le r$.

We summarize the conditions we have imposed so far 
on the parameters $q,r,\qt,\rt$:
\begin{itemize}
\item the sharp case condition $1/q + 1/\qt = 1$;
\item the scaling invariant condition $\beta(q,r;\qt,\rt) = 0$,
  which is equivalent to~\eqref{eq:scaling};
\item conditions for the validity of the local estimates
  with strict inequality,
  \begin{align*}
    & \frac{\sigma - 1}r < \frac\sigma{\rt}, &
    & \frac{\sigma - 1}\rt < \frac\sigma{r}, \\
    & \1q > \sigma \tonde{\1\rt - \1r}, &
    & \1\qt > \sigma \tonde{\1r - \1\rt},
  \end{align*}
  which by~\eqref{eq:scaling} and~\eqref{eq:qqt2} reduce to
  \begin{align*}
    & \frac{(\sigma-1)^2}{\sigma(2\sigma-1)}<\1r<\frac{\sigma-1}{2\sigma-1},&
    & 2\sigma \tonde{\12 - \1r} - 1 < \1q < 2\sigma \tonde{\12 - \1r};
  \end{align*}
\item tecnical conditions needed to recover the $\ell^{r'}$ norm
  from the $\ell^{q'}$ norm, $\qt \le \rt$ and $q \le r$.
\end{itemize}

\section{Applications to {S}chr\"odinger equations}

As an application of theorem~\ref{thm:global}
and theorem~\ref{thm:local}, we derive estimates 
for solutions to inhomogeneous linear Schr\"odinger equations.
Let now $U(t): L^2(\R^n) \to L^2(\R^n)$ be the operator
which describes the solution $u(t,x) = (U(t)f)(x)$
of the homogeneous equation
\begin{equation*}
  \iC \de_t u + \Delta u = 0, \qquad t \in \R, x \in \R^n,
\end{equation*}
with initial data $u(0,x) = f(x)$.
In terms of the Fourier transform we have
\begin{equation}
  \label{eq:uhattxi}
  \mF[U(t)f](\xi) =
  \widehat{u}(t,\xi) = \Eu^{\iC |\xi|^2} \widehat{f}(\xi).
\end{equation}
We also have the explicit formula
\begin{equation}
  \label{eq:utx}
  Tf(t,x) = (U(t)f)(x) = u(t,x) = (4\pi t)^{-n/2}
  \int \Eu^{\iC\abs{x-y}^2/(4t)} f(y) \d y.
\end{equation}
The corresponding $(TT^*)_R$ retarded operator
describes the solution of the inhomogeneous equation
\begin{equation}
  \label{eq:InhomSchrodinger}
  \iC \de_t v + \Delta v = F(t,x), \qquad t>0, x\in \R^n,
\end{equation}
with zero initial data.
We have the explicit formula
\begin{equation}
  \label{eq:vtx}
  (TT^*)_R F(t,x) = v(t,x) = (4\pi)^{-n/2} \int_0^t \int
  \frac{\Eu^{\iC\abs{x-y}^2/(4(t-s))}}{(t-s)^{n/2}}
  F(s,y)\d y\d s.
\end{equation}

Using Plancherel's theorem,
we can immediately verify from~\eqref{eq:uhattxi}
that $U(t)$ satisfies the energy estimate~\eqref{eq:energy}
and also the group property
$U(t) U^*(s) = U(t-s)$.
If we take absolute values
inside the integral in formula~\eqref{eq:utx},
we verify that~\eqref{eq:dispersion} is satisfied 
with $\sigma = n/2$.

Hence, we can apply
the local estimates of theorem~\ref{thm:local}
and the global estimates of theorem~\ref{thm:global}
and obtain the following sufficient conditions.

\begin{proposition} \label{pro:localsuff}
  If $v$ is the solution of~\eqref{eq:InhomSchrodinger}
  with zero initial data
  and inhomogeneous term $F$ supported on $[0,1]\times\R^n$,
  then we have the estimate
  \begin{equation}
    \label{eq:localSchrodinger}
    \norm{v}_{L_t^q([2,3]; L^r(\R^n))} \les
    \norm{F}_{L_t^{\qt'}([0,1]; L^{\rt'}(\R^n))},
  \end{equation}
  whenever $q,r,\qt,\rt$ satisfy the conditions
  \begin{align*}
    & 1 \le q, \qt \le \infty, &
    & 2 \le r, \rt \le \infty, \\
    & \frac{n - 2}r \le \frac{n}{\rt}, &
    & \frac{n - 2}\rt \le \frac{n}{r}, \\
    & \1q \ge \frac{n}{2} \tonde{\1\rt - \1r}, &
    & \1\qt \ge \frac{n}{2} \tonde{\1r - \1\rt},
  \end{align*}
  and if $n = 2$, we must also require $r, \rt < \infty$.
\end{proposition}

\begin{proposition} \label{pro:globalsuff}
  If $v$ is the solution of~\eqref{eq:InhomSchrodinger}
  with zero initial data
  and inhomogeneous term $F$ supported on $\R\times\R^n$,
  then we have the estimate
  \begin{equation}
    \label{eq:globalSchrodinger}
    \norm{v}_{L_t^q(\R; L^r(\R^n))} \les
    \norm{F}_{L_t^{\qt'}(\R; L^{\rt'}(\R^n))},
  \end{equation}
  whenever $(q,r),(\qt,\rt)$ are $n/2$-acceptable pairs
  which satisfy the scaling condition
  \begin{equation*}
    \1q + \1\qt = \frac{n}{2} \tonde{1 - \1r - \1\rt},
  \end{equation*}
  and either the conditions
  \begin{align}
    & \1q + \1\qt < 1, \\
    \label{eq:n2rnrtn2rtnr}
    & \frac{n - 2}r \le \frac{n}{\rt}, \quad
    \frac{n - 2}\rt \le \frac{n}{r},
  \end{align}
  or the conditions
  \begin{align}
    & \1q + \1\qt = 1, \\
    & \frac{n - 2}r < \frac{n}{\rt}, \quad
    \frac{n - 2}\rt < \frac{n}{r}, \\
    & \1r \le \1q, \quad \1\rt \le \1\qt.
  \end{align}
\end{proposition}

We now want to discuss the sharpness of these propositions.
By constructing explicit counterexamples
we are able to show the following necessary conditions.

\begin{proposition} \label{pro:localnecess}
  If the estimate~\eqref{eq:localSchrodinger}
  holds for any $F$ supported on $[0,1]\times\R^n$,
  then $q,r,\qt,\rt$ must satisfy the conditions
  \begin{gather}
    \label{eq:loc-rrt1}
    \1r + \1\rt \le 1, \\
    \label{eq:loc-rrn}
    \abs{\1r - \1\rt} \le \1n, \\
    \label{eq:loc-rqr}
    \frac{n - 2}r - \frac{2}{q} \le \frac{n}{\rt}, \qquad
    \frac{n - 2}\rt - \frac{2}{\qt} \le \frac{n}{r}, \\
    \label{eq:loc-qrr}
    \1q \ge \frac{n}{2} \tonde{\1\rt - \1r}, \qquad
    \1\qt \ge \frac{n}{2} \tonde{\1r - \1\rt}.
  \end{gather}
\end{proposition}

\begin{proposition} \label{pro:globalnecess}
  If the estimate~\eqref{eq:localSchrodinger}
  holds for any $F$ supported on $\R\times\R^n$,
  then $(q,r)$, $(\qt,\rt)$ must be $n/2$-acceptable pairs
  which satisfy the conditions
  \begin{gather}
    \label{eq:glo-scaling}
    \1q + \1\qt = \frac{n}{2} \tonde{1 - \1r - \1\rt}, \\
    \label{eq:glo-qqt1}
    \1q + \1\qt \le 1, \\
    \label{eq:glo-rrn}
    \abs{\1r - \1\rt} \le \1n, \\
    \label{eq:glo-rqr}
    \frac{n - 2}r - \frac{2}{q} \le \frac{n}{\rt}, \qquad
    \frac{n - 2}\rt - \frac{2}{\qt} \le \frac{n}{r}.
  \end{gather}
\end{proposition}

\begin{remark}
  If a bounded linear operator $T:L^p(\R^n) \to L^q(\R^n)$
  is translation invariant,
  \ie\ $T(f \compose \tau_y) = (Tf) \compose \tau_y$
  for any translation $\tau_y(x) = x + y$,
  then we must have that $q$ is bigger or equal to $p$
  (see~\cite{Hor1960}).
  The operator $(TT^*)_R$ defined in~\eqref{eq:vtx}
  has a convolution structure
  and so it is invariant \wrt\ space and time translations.
  As a consequence we obtain the necessity of
  conditions~\eqref{eq:loc-rrt1} and~\eqref{eq:glo-qqt1}.
\end{remark}

\begin{remark}
  The necessity of~\eqref{eq:glo-scaling}
  follows from the scaling properties
  the operator $(TT^*)_R$ defined in~\eqref{eq:vtx}
  under the parabolic scaling
  \begin{equation*}
    (t, x) \leftarrow \TOnde{\frac{t}{\lambda^2}, \frac{x}{\lambda}},
  \end{equation*}
  as $\lambda \to 0^+$ and as $\lambda \to +\infty$.
\end{remark}

Let us now construct some concrete examples of solutions $v(t,x)$
to the inhomogeneous Schr\"odinger equation~\eqref{eq:InhomSchrodinger},
which we use to prove the necessity of the remaining conditions stated
in proposition~\ref{pro:localnecess} and~\ref{pro:globalnecess}.
The \emph{flash} solution of example~\ref{exa:flash}
will correspond to condition~\eqref{eq:loc-qrr};
the \emph{bump} solution of example~\ref{exa:bump}
will prove that $(q,r)$ must be a $n/2$-acceptable pair
in proposition~\ref{pro:globalnecess};
the \emph{focusing} solution of example~\ref{exa:focusing}
will correspond to conditions~\eqref{eq:loc-rqr} and~\eqref{eq:glo-rqr};
the \emph{oscillatory} solution of example~\ref{exa:rrtopt}
will correspond to conditions~\eqref{eq:loc-rrn} and~\eqref{eq:glo-rrn}.

\begin{example} \label{exa:flash}
  Let $\epsi, \eta$ be two small positive parameters
  with $0<\epsi^2<\eta<1$.
  Let $v(t,x)$ be the solution to~\eqref{eq:InhomSchrodinger}
  corresponding to the \emph{flash} forcing term $F$ given by
  the characteristic function 
  \begin{equation*}
    F(s,y) = \chi\tonde{0<s<\epsi^2, \abs{y}<\epsi}.
  \end{equation*}
  When $0<s<\epsi^2, \abs{y}<\epsi$
  and $2<t<3$, $\abs{x} \ll \eta/\epsi$,
  we have
  \begin{equation*}
    \frac{\abs{x-y}^2}{t-s} =
    \tonde{\abs{x}^2 + O(\eta)} \tonde{\1t + O(\epsi^2)} =
    \frac{\abs{x}^2}{t} + O(\eta).
  \end{equation*}
  The oscillating factor in~\eqref{eq:vtx} then becomes
  \begin{equation}
    \label{eq:oscexp}
    \Eu^{\iC\abs{x-y}^2/(4(t-s))} = 
    \Eu^{\iC\abs{x}^2/(4t)} \tonde{1+O(\eta)}.    
  \end{equation}
  Hence, for some small (but fixed) values of $\eta$,
  we can estimate $v$
  in the region $2<t<3$, $\abs{x} \ll \eta/\epsi$,
  \begin{equation*}
    \abs{v(t,x)} \sss
      \int_0^{\epsi^2} \int_{\abs{y}<\epsi} \tonde{1+O(\eta)} \d y \d s \ges
    \epsi^2 \epsi^n.
  \end{equation*}
  We deduce that
  \begin{equation*}
    \frac{\norm{v}_{L_t^q L_x^r}}
    {\norm{F}_{L_s^{\qt'} L_y^{\rt'}}}
    \ges \frac{\epsi^2 \epsi^n \epsi^{-n/r}}{\epsi^{2/\qt'}
      \epsi^{n/\rt'}} =
    \epsi^{2/\qt + n/\rt - n/r}.
  \end{equation*}
  This ratio blows up unless we have 
  \begin{equation*}
    \frac{2}{\qt} \ge n \tonde{\1r - \1\rt},
  \end{equation*}
  which is the necessary condition~\eqref{eq:loc-qrr}
  for the local estimate.
\end{example}

\begin{example} \label{exa:bump}
  Let $\eta$ be a small positive parameter.
  Let $v(t,x)$ be the solution to~\eqref{eq:InhomSchrodinger}
  corresponding to the \emph{bump} forcing term $F$ given by
  the characteristic function 
  \begin{equation*}
    F(s,y) = \chi\tonde{0<s<1, \abs{y}<1}.
  \end{equation*}
  When $0<s<1, \abs{y}<1$
  and $t>2$, $\abs{x} \ll \eta t$,
  we have
  \begin{equation*}
    \frac{\abs{x-y}^2}{t-s} =
    \tonde{\abs{x}^2 + O(\eta t)} \tonde{\1t + O\tonde{\1{t^2}}} =
    \frac{\abs{x}^2}{t} + O(\eta).
  \end{equation*}
  The oscillating factor in~\eqref{eq:vtx}
  still behaves as in~\eqref{eq:oscexp}.
  Hence, for some small (but fixed) values of $\eta$,
  we can estimate $v$
  in the region $t>2$, $\abs{x} \ll \eta t$,
  \begin{equation*}
    \abs{v(t,x)} \sss t^{-n/2} 
      \int_0^{1} \int_{\abs{y}<1} \tonde{1+O(\eta)} \d y \d s \ges t^{-n/2}.
  \end{equation*}
  We deduce that
  \begin{equation*}
    \norm{v(t)}_{L^r} \ges t^{-n/2} t^{n/r} = t^{-n(1/2 - 1/r)}.
  \end{equation*}
  When $v \in \LtLX{q}{r}$ then
  the right hand side must belong to $L^q(\R)$,
  but this happens only if $(q,r) = (\infty, 2)$
  or if the integrability condition
  \begin{equation*}
    -q n \tonde{\12 - \1r} < -1,
  \end{equation*}
  is satisfied,
  which is equivalent to say that $(q,r)$ is a $n/2$-acceptable pair.
\end{example}

\begin{example} \label{exa:focusing}
  Let $\epsi, \eta$ be two small positive parameters
  with $0<\epsi^2<\eta<1$.
  Let $v(t,x)$ be the solution to~\eqref{eq:InhomSchrodinger}
  corresponding to the \emph{focusing} forcing term $F$
  given by the characteristic function 
  \begin{equation*}
    F(s,y) =
    \chi\tonde{0<s<\epsi^2,
      \abs{\abs{y}-\frac\eta\epsi}<\epsi}.
  \end{equation*}
  When $0<s<\epsi^2, \abs{\abs{y}-\frac\eta\epsi}<\epsi$
  and $2<t<3$, $\abs{x} \ll \epsi$,
  we have
  \begin{equation*}
    \frac{\abs{x-y}^2}{t-s} =
    \tonde{\frac{\eta^2}{\epsi^2} + O(\eta)}
    \tonde{\1t + O(\epsi^2)}
    = \frac{\eta^2}{\epsi^2 t} + O(\eta).
  \end{equation*}
  The oscillating factor in~\eqref{eq:vtx} then becomes
  \begin{equation*}
    \Eu^{\iC\abs{x-y}^2/(4(t-s))} =
    \Eu^{\iC\eta^2/(4\epsi^2t)} \tonde{1+O(\eta)}.    
  \end{equation*}
  Hence, for some small (but fixed) values of $\eta$,
  we can estimate $v$
  in the region $2<t<3$, $\abs{x} \ll \epsi$,
  \begin{equation*}
    \abs{v(t,x)} =
    \int_0^{\epsi^2} \int_{\abs{\abs{y}-\frac\eta\epsi}<\epsi}
    \tonde{1+O(\eta)} \d y \d s \ges
    \epsi^2 \epsi^{-n+2}.
  \end{equation*}
  We deduce that
  \begin{equation*}
    \frac{\norm{v}_{L_t^q L_x^r}}
    {\norm{F}_{L_s^{\qt'} L_y^{\rt'}}}
    \ges \frac{\epsi^2 \epsi^{-n+2} \epsi^{n/r}}
    {\epsi^{2/\qt'} \epsi^{(-n+2)/\rt'}} =
    \epsi^{n/r + 2/\qt + (n-2)/\rt}.
  \end{equation*}
  This ratio blows up unless we have 
  \begin{equation*}
    \frac{n-2}{\rt} - \frac2\qt \le \frac{n}{r},
  \end{equation*}
  which is the necessary condition~\eqref{eq:loc-rqr}
  for the local estimate
  or the necessary condition~\eqref{eq:glo-rqr} 
  for the global estimate.
\end{example}

\begin{example} [\cite{Fos2003tri}] \label{exa:rrtopt}
  Let $R \gg 1$ and $0< \eta <1$.
  We choose
  \begin{equation*}
    F(s,y) = \Eu^{\iC 2 R^2 s^2}
    \chi\tonde{0<s<1, \abs{y} \le \frac\eta{R}}.
  \end{equation*}
  We have 
  $\norm{F}_{L^{\rt'}(\R^n; L^{\qt'}([0,1]))} \sss
  (\eta/R)^n$.
  We can write the solution $v$ as
  \begin{equation}
    \label{eq:38}
    v(t,x) = \int_{|y|<\eta/R} \int_0^1
    \frac{\Eu^{\iC \tonde{2 R^2 s^2 - |x-y|^2/(4(t-s))}}}
    {(t-s)^{n/2}} \d s \d y
    = \int_{|y|<\eta/R} I\TOnde{t, \frac{x-y}{2R}, R} \d y,
  \end{equation}
  where $I$ is the oscillatory integral
  \begin{equation*}
    I(t,z,R) =
    \int_0^1 \Eu^{\iC R^2 \varphi(s;t,z)} \psi(s;t) \d s,
  \end{equation*}
  with phase $\varphi(s;t,z) = 2s^2 - |z|^2/(t-s)$
  and amplitude $\psi(s;t) = 1/(t-s)^{n/2}$.
  The first and second derivatives of the phase are
  \begin{align*}
    \de_s \varphi &= 4s - \frac{|z|^2}{(t-s)^2}, &
    \de_s^2 \varphi &= 4 - \frac{2|z|^2}{(t-s)^3}.
  \end{align*}
  When $t \in [2,3]$ and $1/2 \le |z| \le 1$,
  derivatives with respect to $s$ of all orders
  for $\varphi$ and $\psi$
  are uniformly bounded by absolute constants
  and the phase $\varphi$
  has exactly one non degenerate critical point
  $s_* = s_*(t,z)$ in $[0,1]$,
  \begin{align*}
    s_* &= \frac{|z|^2}{4(t-s_*)^2} \in \quadre{\1{48},\14}, &
    \de_s \varphi(s_*, t, z) &= 0, &
    \de_s^2 \varphi(s_*, t, z) &= 4 - \frac{2s_*}{t-s_*}
    \in \quadre{\frac23,12}.
  \end{align*}
  By standard stationary phase methods
  (see~\cite{Ste1993, Fos2003tri}),
  we obtain that the integral $I$
  decays like $1/R$ as $R \to \infty$,
  more precisely
  \begin{equation*}
    I(t,z,R) =
    \frac{J_*(t,z) \Eu^{\iC R^2 \varphi_*(t,z)}}{R}
    + O\tonde{\1{R^2}},
  \end{equation*}
  where
  \begin{align*}
    J_*(t,z) &= \Eu^{i\pi/4} \psi(s_*(t,z);t)
    \sqrt{\frac{2\pi}{\de_s^2 \varphi(s_*(t,z);t,z)}}, &
    \varphi_*(t,z) &= \varphi(s_*(t,z);t,z).
  \end{align*}
  By the above computations, when $t \in [2,3]$ and $|z|<1$
  we have
  \begin{equation}
    \label{eq:39}
    \abs{J_*(t,z)} \sss 1;
  \end{equation}
  moreover, 
  \begin{equation*}
    \nabla_z \varphi_*(t,z) = \nabla_z \varphi(s_*;t,z) =
    - \frac{2z}{t-s_*} = O(1),
  \end{equation*}
  so that we have
  $\varphi_*(t,z) - \varphi(t,z_0) = O(z-z_0)$.
  In particular, this shows that the oscillatory factor
  \begin{equation*}
     \Eu^{\iC R^2 \varphi_*\tonde{t, \frac{x-y}{2R}}}
     = \Eu^{\iC R^2 \varphi_*\tonde{t,\frac{x}{2R}} + O(\eta)}
     = \Eu^{\iC R^2 \varphi_*\tonde{t,\frac{x}{2R}}}
     \tonde{1+O(\eta)},
  \end{equation*}
  does not oscillates too much 
  when $|y| \le \eta/R$, $R \le |x-y| \le 2R$
  and $\eta$ is sufficiently small.
  It follows that
  \begin{equation*}
    I\tonde{t, \frac{x-y}{2R}, R} = 
    \frac{\Eu^{\iC R^2 \varphi_*(t,x/(2R))}}{R}
    J_*\tonde{t, \frac{x-y}{2R}} \tonde{1 + O(\eta)}
    + O\tonde{\1{R^2}},
  \end{equation*}
  which inserted in~\eqref{eq:38} and using~\eqref{eq:39}
  proves that
  \begin{equation*}
    \abs{v(t,x)} \ges \frac{\eta^n}{R^{n+1}}
  \end{equation*}
  on the region $2\le t\le3$, $R+\eta/R<|x|<2R-\eta/R$.
  Thus, the ratio
  \begin{equation*}
    \frac{\norm{v}_{L^r(\R^n; L^q([2,3]))}}
    {\norm{F}_{L^{\rt'}(\R^n; L^{\qt'}([0,1]))}} \ges
    \frac{R^{-n-1} \cdot R^{n/r}}{R^{-n/\rt'}} = 
    R^{n(1/r - 1/\rt) - 1},
  \end{equation*}
  cannot be bounded as $R \to \infty$ unless
  \begin{equation*}
    \1r - \1\rt \le \1n,
  \end{equation*}
  which is the necessary condition~\eqref{eq:loc-rrn}
  for the local estimate
  or the necessary condition~\eqref{eq:glo-rrn}
  for the global estimate.
\end{example}

\section{Open questions}

\begin{remark}
  For the local estimates~\eqref{eq:localSchrodinger},
  the range of values for the exponents $(q,r;\qt,\rt)$
  which is described by the necessary conditions
  of proposition~\ref{pro:localnecess}
  is larger than the corresponding range described
  by the sufficient conditions
  of proposition~\ref{pro:localsuff}.
  In particular, our examples do not exclude the possibility
  that local estimates could be valid for some $(q,\qt)$
  when $(r,\rt)$ are in one of the following ranges, $R_j$,
  which are not covered by proposition~\ref{pro:localsuff}:
  \begin{align*}
    & R_1: \quad
    \1r > \12, \quad \1r + \1\rt \le 1, \quad \1r - \1\rt \le \1n; \\
    & R_2: \quad
    \1\rt > \12, \quad \1r + \1\rt \le 1, \quad \1\rt - \1r \le \1n; \\
    & R_3: \quad \frac{n-2}{r} > \frac{n}{\rt}, \quad \1r - \1\rt \le \1n; \\
    & R_4: \quad \frac{n-2}{\rt} > \frac{n}{r}, \quad \1\rt - \1r \le \1n.
  \end{align*}
  To our knowledge, this is still an open problem.
  It seems that using only the conservation of energy
  and the dispersive properties of the linear evolution
  for the homogeneous Schr\"odinger equation,
  as we do with Strichartz estimates,
  is not enough to reach these regions.
  Is there some other smoothing property associated with
  the inhomogeneous operator~\eqref{eq:vtx}?
  
  If, instead of using $\LtLX{q}{r}$ norms,
  we measure functions using $L^r(X;L^q(\R))$ norms,
  then we are able to prove estimates 
  for $(r,\rt)$ inside the regions $R_3$ and $R_4$;
  details on this can be found in~\cite{Fos2003tri}.
\end{remark}

\begin{figure}
  \centering
  \begin{picture}(0,0)%
    \includegraphics{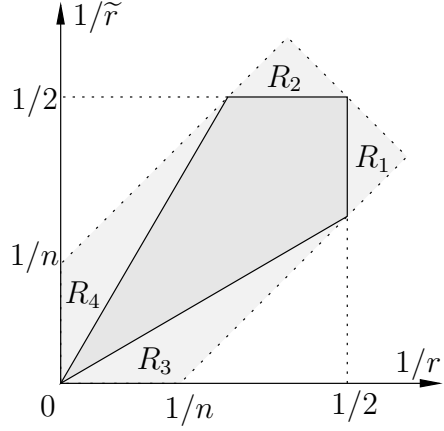}%
  \end{picture}%
  \begin{picture}(2742,2680)(43,-1979)
    \put(1999,-1902){\makebox(0,0)[lb]{\smash{{\SetFigFont{12}{14.4}{%
              \familydefault}{\mddefault}{\updefault}{$1/2$}}}}}
    \put( 952,-1930){\makebox(0,0)[lb]{\smash{{\SetFigFont{12}{14.4}{%
              \familydefault}{\mddefault}{\updefault}{$1/n$}}}}}
    \put( -20, -968){\makebox(0,0)[lb]{\smash{{\SetFigFont{12}{14.4}{%
              \familydefault}{\mddefault}{\updefault}{$1/n$}}}}}
    \put(2146, -351){\makebox(0,0)[lb]{\smash{{\SetFigFont{12}{14.4}{%
              \familydefault}{\mddefault}{\updefault}{$R_1$}}}}}
    \put(1592,  156){\makebox(0,0)[lb]{\smash{{\SetFigFont{12}{14.4}{%
              \familydefault}{\mddefault}{\updefault}{$R_2$}}}}}
    \put( 779,-1618){\makebox(0,0)[lb]{\smash{{\SetFigFont{12}{14.4}{%
              \familydefault}{\mddefault}{\updefault}{$R_3$}}}}}
    \put( 326,-1204){\makebox(0,0)[lb]{\smash{{\SetFigFont{12}{14.4}{%
              \familydefault}{\mddefault}{\updefault}{$R_4$}}}}}
    \put( -10,   17){\makebox(0,0)[lb]{\smash{{\SetFigFont{12}{14.4}{%
              \familydefault}{\mddefault}{\updefault}{$1/2$}}}}}
    \put( 172,-1915){\makebox(0,0)[lb]{\smash{{\SetFigFont{12}{14.4}{%
              \familydefault}{\mddefault}{\updefault}{$0$}}}}}
    \put(2401,-1636){\makebox(0,0)[lb]{\smash{{\SetFigFont{12}{14.4}{%
              \familydefault}{\mddefault}{\updefault}{$1/r$}}}}}
    \put( 376,  539){\makebox(0,0)[lb]{\smash{{\SetFigFont{12}{14.4}{%
              \familydefault}{\mddefault}{\updefault}{$1/\rt$}}}}}
  \end{picture}%
  \caption{Necessary and sufficient conditions on $r$ and $\rt$
    for local inhomogeneous estimates for the Schr\"odinger equation.}
  \label{fig:Srtr}
\end{figure}

\begin{remark}
  Similarly,
  for the global estimates \eqref{eq:globalSchrodinger},
  the range of values for the exponents $(q,r;\qt,\rt)$
  which is described by the necessary conditions
  of proposition~\ref{pro:globalnecess}
  is larger than the corresponding range described
  by the sufficient conditions
  of proposition~\ref{pro:globalsuff}.
  The gap here lies in the difference between 
  condition~\eqref{eq:n2rnrtn2rtnr} and condition~\eqref{eq:glo-rqr}.
\end{remark}

\begin{remark}
  Finally, one last question.
  If the inhomogeneous Schr\"odinger equation turns out to have
  better integrability properties than the ones provided by 
  proposition~\ref{pro:localsuff} or proposition~\ref{pro:globalsuff},
  is it then possible to construct some family $U(t)$ of evolution operators
  which satisfy the requirements~\eqref{eq:energy} and~\eqref{eq:dispersion}
  and such that their integrability properties
  for the inhomogeneous estimates
  are \emph{exactly} those given by 
  theorem~\ref{thm:local} and theorem~\ref{thm:global}?
\end{remark}

\appendix

\section{Details on the set $\mE_*$
  defined in section~\ref{sec:proof-local-estim}}
\label{sec:details-about-set}

We construct the set $\mE_*$ in three steps.

First, $\mE_*$ contains the point $(1/Q,1/R; 1/\Qt, 1/\Rt)$
when the pairs $(Q,R)$ and $(\Qt, \Rt)$ are sharp admissible.
This is a square in $[0,1]^4$ defined by the equations
\begin{align}
  \label{eq:QRQR}
  &\1Q = \sigma \tonde{\12 - \1R}, &
  &\1\Qt = \sigma \tonde{\12 - \1\Rt}, &
  & 0 \le \1Q, \1\Qt, \1R, \1\Rt \le \12,
\end{align}
and, if $\sigma=1$, we must also require
$(Q,R) \ne (2,\infty)$ and $(\Qt, \Rt) \ne (2,\infty)$.

Second, $\mE_*$ contains the convex hull
of the above square with the point
$(1/\infty, 1/\infty; 1/\infty, 1/\infty)$.
These are points of the form
$(\theta/Q, \theta/R; \theta/\Qt, \theta/\Rt)$
where $(Q,R)$ and $(\Qt, \Rt)$ satisfy~\eqref{eq:QRQR}
and $0 \le \theta \le 1$.

Third, $\mE_*$ contains points of the form
$(1/q, 1/r; 1/\qt, 1/\rt)$ where
\begin{align*}
  &\1q \ge \frac\theta{Q}, &
  &\1r = \frac\theta{R}, & 
  &\1\qt \ge \frac\theta\Qt, &
  &\1\rt = \frac\theta\Rt.
\end{align*}

Hence, the set $\mE_*$ is the set of points
$(1/q, 1/r; 1/\qt, 1/\rt) \in [0,1]^4$
for which there exist $Q,R,\Qt,\Rt,\theta$ such that
\begin{align*}
  &\1Q = \sigma \tonde{\12 - \1R}, &
  &\1\Qt = \sigma \tonde{\12 - \1\Rt}, \\
  & 0 \le \1Q, \1\Qt \le \12, &
  & 0 \le \1R, \1\Rt \le \12, &
  & 0 \le \theta \le 1, \\
  &\1q \ge \frac\theta{Q}, &
  &\1\qt \ge \frac\theta\Qt, &
  &\1r = \frac\theta{R}, & 
  &\1\rt = \frac\theta\Rt,
\end{align*}
and if $\sigma=1$ we must also require
that $R \ne \infty$ and $\Rt \ne \infty$.

Using the last two equalities,
we can eliminate $R$ and $\Rt$,
\begin{align*}
  &\frac\theta{Q} = \sigma \tonde{\frac\theta2 - \1r}, &
  &\frac\theta\Qt = \sigma \tonde{\frac\theta2 - \1\rt}, \\
  & 0 \le \frac\theta{Q}, \frac\theta\Qt \le \frac\theta2, &
  & 0 \le \1r, \1\rt \le \frac\theta2, &
  & 0 \le \theta \le 1, \\
  &\1q \ge \frac\theta{Q}, &
  &\1\qt \ge \frac\theta\Qt.
\end{align*}
Using the first two equalities,
we can eliminate $Q$ and $\Qt$,
\begin{align*}
  & 0 \le \sigma \tonde{\frac\theta2-\1r} \le \frac\theta2, &
  & 0 \le \sigma \tonde{\frac\theta2-\1\rt} \le
  \frac\theta2, \\
  & 0 \le \1r, \1\rt \le \frac\theta2, &
  & 0 \le \theta \le 1, \\
  &\1q \ge \sigma \tonde{\frac\theta2 - \1r}, &
  &\1\qt \ge \sigma \tonde{\frac\theta2 - \1\rt}.
\end{align*}
We rearrange these inequalities,
\begin{align*}
  & 0 \le \theta \le 1, \\
  & \frac{\sigma-1}\sigma \frac\theta2 \le \1r \le
  \frac\theta2, &
  & \frac{\sigma-1}\sigma \frac\theta2 \le \1\rt \le
  \frac\theta2, \\
  & \frac{\sigma\theta}2 \le \1q + \frac\sigma{r}, &
  & \frac{\sigma\theta}2 \le \1\qt + \frac\sigma\rt.
\end{align*}
We isolate the quantity $1/\theta$,
\begin{align*}
  1 \le &\1\theta \le \infty, \\
  \frac{\sigma-1}\sigma \frac{r}2 \le &\1\theta \le
  \frac{r}2, &
  \frac{\sigma-1}\sigma \frac\rt2 \le &\1\theta \le
  \frac\rt2, \\
  \frac{\frac\sigma2}{\1q + \frac\sigma{r}} \le &\1\theta, &
  \frac{\frac\sigma2}{\1\qt + \frac\sigma\rt} \le &\1\theta.
\end{align*}
There exists some $\theta$ which satisfies 
the above system of inequalities if and only if 
each expression on the left of $1/\theta$
is less or equal to each expression on the right.
This means that we must have
\begin{align*}
  1 &\le \frac{r}2, &
  1 &\le \frac\rt2, \\
  \frac{\sigma-1}\sigma \frac{r}2 &\le \frac\rt2, &
  \frac{\sigma-1}\sigma \frac\rt2 &\le \frac{r}2, \\
  \frac\sigma2 &\le \frac{r}2 \tonde{\1\qt+\frac\sigma\rt}, &
  \frac\sigma2 &\le \frac\rt2 \tonde{\1q+\frac\sigma{r}}.
\end{align*}
We rearrange these inequalities in a final form
\begin{align*}
  r &\ge 2, &
  \rt &\ge 2, \\
  \frac{\sigma-1}\rt &\le \frac\sigma{r}, &
  \frac{\sigma-1}r &\le \frac\sigma\rt, \\
  \frac\sigma{r} - \frac\sigma\rt &\le \1\qt, &
  \frac\sigma\rt - \frac\sigma{r} &\le \1q.
\end{align*}
We should also remember that in the case $\sigma=1$
we had to exclude the case $r=\infty$ or $\rt=\infty$.
These are the conditions which describe the set $\mE_*$ 
and are the same which appear in the statement of theorem~\ref{thm:local}.


\begin{thebibliography}{1}
  
\bibitem{BerLof1976}
  J{\"o}ran Bergh and J{\"o}rgen L{\"o}fstr{\"o}m,
  \emph{Interpolation spaces. {A}n introduction},
  Springer-Verlag, Berlin, 1976,
  Grundlehren der Mathematischen Wissenschaften, No. 223.
  
\bibitem{Fos2003tri}
  Damiano Foschi,
  \emph{Some remarks on the $l^p-l^q$ boundedness
    of trigonometric sums and oscillatory integrals},
  preprint 2003.
  
\bibitem{Har1990}
  Joergen Harmse,
  \emph{On {L}ebesgue space estimates for the wave equation},
  Indiana Univ. Math. J. \textbf{39} (1990), no.~1, 229--248.
  
\bibitem{Hor1960}
  Lars H{\"o}rmander,
  \emph{Estimates for translation invariant operators in {$L^p$} spaces},
  Acta Math. \textbf{104} (1960), 93--140.
  
\bibitem{Kat1994}
  Tosio Kato,
  \emph{An {$L\sp {q,r}$}-theory for nonlinear {S}chr\"odinger equations},
  Spectral and scattering theory and applications,
  Adv. Stud. Pure Math., vol.~23, Math. Soc. Japan, Tokyo, 1994, pp.~223--238.
  
\bibitem{KeeTao1998}
  Markus Keel and Terence Tao,
  \emph{Endpoint {S}trichartz estimates},
  Amer. J. Math. \textbf{120} (1998), no.~5, 955--980.
  
\bibitem{Obe1989}
  Daniel~M. Oberlin,
  \emph{Convolution estimates for some distributions
    with singularities on the light cone},
  Duke Mathematical Journal \textbf{59} (1989), no.~3, 747--757.
  
\bibitem{Ste1970}
  Elias~M. Stein,
  \emph{Singular integrals and differentiability properties of functions},
  Princeton University Press, 1970.
  
\bibitem{Ste1993}
  Elias~M. Stein,
  \emph{Harmonic analysis:
    real-variable methods, orthogonality, and oscillatory integrals},
  Princeton University Press, Princeton, NJ, 1993.
  
\end{thebibliography}
\end{document}